\documentclass[12pt,draft]{amsart}
\usepackage{amsfonts,amscd,latexsym}
\title[Birational properties...]
{Birational properties of pencils of del Pezzo surfaces of degree 1 and 2.}
\thanks{This work was partially supported by the grants RFBR no. 99--01--01132,
Grant of Leading Scientific Schools no. 96--15--96146, and INTAS-OPEN 97/2072}
\author{Mikhail Grinenko}
\address{Steklov Mathematical Institute}
\address{Max-Planck Institut f\"ur Mathematik}
\email{grin@mi.ras.ru / grinenko@mpim-bonn.mpg.de}
\date{}

\newtheorem{theorem}{\sc Theorem}[section]
\newtheorem{proposition}[theorem]{\sc Proposition}
\newtheorem{lemma}[theorem]{\sc Lemma}
\newtheorem{corollary}[theorem]{\sc Corollary}

\newtheorem{conjecture}[theorem]{\sc Conjecture}
\newtheorem{definition}[theorem]{\sc Definition}

\makeatletter

\@addtoreset{equation}{section}
\newcommand{\l@abcd}[2]{\hbox to\textwidth{#1\dotfill #2}}

\makeatother

\newcommand*{\mybegintheorem}[1]{\begin{trivlist}\it%
      \item[\hspace{\labelsep}{\bf #1}]}
\newcommand*{\myendtheorem}{\end{trivlist}}
\newenvironment*{theorem*}{\mybegintheorem{Theorem.}}{\myendtheorem}
\newenvironment*{proposition*}{\mybegintheorem{Proposition.}}{\myendtheorem}
\newenvironment*{corollary*}{\mybegintheorem{Corollary.}}{\myendtheorem}
\newenvironment*{definition*}{\mybegintheorem{Definition.}}{\myendtheorem}

\theoremstyle{remark}
\newtheorem{remark}[theorem]{\sc Remark}

\renewcommand{\phi}{\varphi}
\renewcommand{\epsilon}{\varepsilon}

\newcommand{\PS}{{{\mathbb P}^6}}
\newcommand{\PF}{{{\mathbb P}^5}}
\newcommand{\PQ}{{{\mathbb P}^4}}
\newcommand{\PT}{{{\mathbb P}^3}}
\newcommand{\PTw}{{{\mathbb P}^2}}
\newcommand{\POn}{{{\mathbb P}^1}}
\newcommand{\ZA}{{\mathbb Z}}
\newcommand{\QA}{{\mathbb Q}}
\newcommand{\RA}{{\mathbb R}}

\newcommand{\mc}{\mathcal}

\newcommand{\eqdef}{\stackrel{\rm def}{=}}

\newcommand{\Supp}{\mathop{\rm Supp}\nolimits}
\newcommand{\mult}{\mathop{\rm mult}\nolimits}

\newcommand{\Pic}{\mathop{\rm Pic}\nolimits}

\newcommand{\Bas}{\mathop{\rm Bas}\nolimits}

\newcommand{\Proj}{\mathop{{\bf Proj}\:}\nolimits}

\newcommand{\eps}{\varepsilon}

\begin{document}
\abstract{In this paper we study the birational rigidity problem for 
smooth Mori fibrations on del Pezzo surfaces of degree 1 and 2. For
degree 1 we obtain a complete description of rigid and non-rigid cases.}
\endabstract
\maketitle

\section{Notions and results.}
\label{sec1}

Due to progress in the Minimal Model Program, the classification
problem in the modern birational geometry for varieties of negative
Kodaira dimension can be formulated as follows: given a class of birational
equivalency, to describe all Mori fibrations and birational maps between
them.

Recall that a normal variety $V$ with at most $\QA$-factorial terminal
singularities is called {\it a Mori fibration} if there exists an extremal
contraction $\phi: V\to S$ of fibering type, i.e.

\begin{itemize}
\item[1)] $\phi$ is a morphism with connected fibers onto a normal variety
$S$ and $\dim S < \dim V$;
\item[2)] $-K_V$ is $\phi$-ample and the relative Picard number is equal to 1:
$\rho(V/S)=\rho(V)-\rho(S)=1$.
\end{itemize}

Note that often comparing birational classes of varieties, we suffice to know
Mori structures rather than Mori fibrations themselves (roughly speaking,
Mori fibrations modulo birational maps over the base). From this viewpoint,
the following class of varieties is rather important (and simple for 
describing):

\begin{definition}
A Mori fibration $V/S$ is said to be birationally rigid, if any birational
map $\chi: V\dasharrow V'$ onto another Mori fibration $V'/S'$ is birational
over the base ({\it "square"}), i.e., there exists a birational map
$\psi: S\dasharrow S'$ making the following diagram to be commutative:
$$
\begin{array}{ccc}
 V &\stackrel{\chi}{\dasharrow} & V' \\
 \downarrow && \downarrow \\
 S & \stackrel{\psi}{\dasharrow} & S'
\end{array}
$$
\end{definition}

Now let $\rho: V\to S$ be a Mori fibration. Then any non-empty linear system
${\mc D}$ on $V$ is a subsystem of $|a(-K_V)+\rho^*(A)|$, where 
$a=\mu(V,{\mc D})\ge 0$ is the so-called {\it quasi-effective threshold}, and
$A$ is a divisor on $S$. The behaviour of linear systems under birational
maps between Mori fibrations is given by the following theorem (\cite{Corti1},
theorem 4.2):

\begin{proposition}
\label{th_Corti}
Let $\rho: V\to S$ and $\rho': V'\to S'$ be Mori fibrations, 
$\chi: V\dasharrow V'$ a birational map, ${\mc D'}=|n'(-K_{V'})+\rho'{}^*(A')|$
a very ample linear system on $V'$, where $A'$ is an ample divisor on $S'$.
Denote ${\mc D}=\chi_*^{-1}{\mc D'}\subset|n(-K_V)+\rho^*(A)|$. Then:
\begin{itemize}
\item[(i)] $n\ge n'$ always, and if $n=n'$ then $\chi$ is square;
\item[(ii)] if a log pair $K_V+\frac1n{\mc D}$ is canonical and numerically
effective, then $\chi$ is an isomorphism inducing also an isomorphism between
the bases.
\end{itemize}
\end{proposition}

The statement ${\it (i)}$ of the proposition allows us to describe
various Mori structures on $V$, while ${\it (ii)}$ can be used to get
Mori models (i.e., various Mori fibrations birational to $V$).

There is another important number arising from the maximal singularities
method (now playing the key role in the birational classification). This is
so-called {\it the adjunction threshold}.

\begin{definition}
A non-negative number $\alpha(V,{\mc D})$ is called the adjunction threshold of
a pair $(V,{}\mc D)$, where ${\mc D}$ is a non-empty linear system without
fixed components, $V$ is projective and non-singular in codimension 1,
if $\alpha(V,{\mc D})$ is the smallest number such that for any positive 
integers $m$ and $n$ with $\frac{m}{n}>\alpha(V,{\mc D})$ the linear system
$|nD+mK_V|$ is empty, where $D$ is an element of ${\mc D}$.
\end{definition}

The following statement clarifies the role of the adjunction threshold 
(\cite{Pukh2}, 2.1.):

\begin{proposition}
Let $\chi: V\to V'$ be a birational map between terminal $\QA$-factorial
varieties,  ${\mc D'}$ a non-empty linear system on $V'$ without fixed
components, and ${\mc D}=\chi_*^{-1}{\mc D'}$ its strict transform.
If $\alpha(V,{\mc D})>\alpha(V',{\mc D'})$, then a pair 
$K_V+\frac1{\alpha(V,{\mc D})}{\mc D}$ is not canonical (in other words,
${\mc D}$ has a maximal singularity)
\end{proposition}

This condition "to be not canonical" is strong enough, and we can often
either prove its impossibility, or decrease the adjunction threshold 
applying an "untwisting" birational automorphism or jumping to another
Mori model (as in the Sarkisov program). A relation between 
both the thresholds is given by the following obvious lemma:

\begin{lemma}
\label{easy_lemma}
Let ${\mc D}\subset|n(-K_V)+\rho^*(A)|$ be a non-empty linear system
without fixed components on a Mori fibration $\rho:V\to S$. Always
$n=\mu(V,{\mc D})\ge \alpha(V,{\mc D})$, and if $A$ is effective, 
then the equality holds; if $C\circ\rho^*(A)<0$ for a general curve $C$
covering the base, then $\mu(V,{\mc D})>\alpha(V,{\mc D})$.
\end{lemma}

The key idea is that the relation between these thresholds determines
whether a variety is rigid. For Mori fibrations that are pencils of del 
Pezzo surfaces, it looks as follows:

\begin{conjecture}
\label{main_conj}
Let $V/\POn$ be a pencil of del Pezzo surfaces of degree 1,2 or 3 
(we assume it to be a Mori fibration). $V/\POn$ is birationally rigid
if and only if $\alpha(V,{\mc D})=\mu(V,{\mc D})$ for any linear system
without fixed components (i.e., a linear system $|n(-K_V)-F|$ either is
empty, or has a base divisor; here $F$ is the class of a fiber).
\end{conjecture}

Let us note that "nearly all" (smooth) pencils of del Pezzos of degree 1,
2 or 3 are sufficiently "twisted" along the base, and in \cite{Pukh1}
their rigidity is proved. In this paper, we deal with all cases of degree 1
(theorem \ref{dp1_th}) and advance in degree 2 (theorem \ref{dp2_th}). 
In particular, conjecture \ref{main_conj} turns out to be true for smooth 
cases of degree 1.

The paper is based on using the maximal singularities method (\cite{IskMan},
\cite{Isk2}) in its most perfect kind \cite{Pukh2}.

Everywhere in this paper the characteristic of the ground field is assumed
to be 0.

The author would like to thank the Directors and the stuff of the Max-Planck
Institute for Mathematics in Bonn for hospitality and excellent conditions
during the work.

\section{Smooth varieties with a pencil of del Pezzo surfaces of degree 1.}
\label{sec2}

\subsection{The essential construction.}
\label{subsec1}
We follow to the idea given in \cite{Isk1}.

Let $\rho: V\to \POn$ be a smooth Mori fibration on del Pezzo surfaces of
degree 1, $\Pic V = \ZA[-K_V]\oplus\ZA[F]$, where $F$ is the class of a 
fiber. By the Grauert theorem, $\rho_*{\mc O}(-2K_V+mF)$ is a vector bundle
of rank 4 on $\POn$. We may choose $m$ such that
$$
  {\mc E}=\rho_*{\mc O}(-2K_V+mF)\simeq {\mc O}\oplus{\mc O}(n_1)
    \oplus{\mc O}(n_2)\oplus{\mc O}(n_3),
$$
where $0\le n_1\le n_2\le n_3$. Suppose $b=n_1+n_2+n_3$.

Note that there exists a special section $s_B$:
$$
       s_B=Bas|-K_V+kF|
$$
for $k\gg 0$. Obviously, any fiber $S\in|F|$ is smooth at the point 
$s_B\cap S$, which is also a unique base point of $|-K_S|$.

Let $\pi: X\to\POn$ be a natural projection for $X\eqdef\Proj{\mc E}$,
$\phi: V\to X$ a morphism defined by $|-2K_V+mF|$. 
Then $\rho=\pi\circ\phi$. It is easy to see that if $Q\subset X$ is the 
image of $\phi$, then $\phi: V\to Q$ is a degree 2 finite morphism branched
along a smooth divisor $R_Q=R\cap Q$ that does not intersect 
$t_B\eqdef\phi(s_B)$, where $R\subset X$ is a divisor.

$Q$ is restricted to any fiber $T$ of $X/\POn$ as a non-degenerated quadratic
cone with $T\cap t_B$ as its vertex, and $R|_T$  is a cubic surface.
So if $M$ is the class of the tautological bundle on $X$ and $L$ is the class
of a fiber, we have
$$
\begin{array}{cc}
   Q=2M+aL;  &  R=3M+cL.
\end{array}
$$

Let $t_0$ be a section corresponding to the surjection 
${\mc E}\to{\mc O}\to 0$, $l$  the class of a line in $L\cong\PT$; 
set $t_B\sim t_0+\eps l$, where $\eps$ is obviously nonnegative.

Choose a very ample divisor $H\sim M+\beta l$, $\beta>0$. So $D=H\cap Q$
is a conic bundle with $N=t_B\circ H$ degenerations. We may suppose that
all fibers of D are reduced and all its singular points are du Val 
${\bf A_1}$.
It is easy to compute that $K_D^2=8-2b-2\beta-3a$. Let us blow up all 
singular points of $D$ and then contract all (-1)-curves (there will be 
$2N$ of them) onto a ruled surface $D'$. We have
$$
     K_D^2+2N=K_{D'}^2=8,
$$
hence
$$
     N=b+\beta+\frac32a.
$$
Note that $a$ must be even. Set $a=2a'$. Then,
$$
   N=t_{б}\circ H=(t_0+\eps l)\circ (M+\beta l)=\beta+\eps,
$$
and we obtain $\eps=b+3a'$. Since $R|_Q\cap t_B=\emptyset$, we have
$0=(t_0+\eps l)\circ(3M+cL)=3\eps+c$. So, we have proved the following
lemma:

\begin{lemma}
The following equalities hold:
$$
\begin{array}{l}
Q=2M-\frac23(b-\eps)L, \\
t_B=t_0+\eps l,  \\
R=3M-3\eps L.
\end{array}
$$
\end{lemma}

In the sequel we will always assume $b>0$. Really, if $b=0$, then
$Q=2M$ and $\eps=0$, whence $V\cong\POn\times S$ for some smooth
del Pezzo surface $S$. So $V$ is rational and even not a Mori fibration
because of the Picard group rank.

\begin{lemma}
\label{struct_lemma}
The only following two cases may occur:
\begin{itemize}
\item[1)] $\eps=0$, $2n_2=n_1+n_3$, $n_1$ and $n_3$ are even; so we have
$$
\begin{array}{ccc}
Q=2M-2n_2L, & t_B=t_0, & R=3M;
\end{array}
$$
\item[2)] $\eps=n_1>0$, $n_3=2n_2$, $n_1$ is even and $n_2\ge 3n_1$;
in this case
$$
\begin{array}{ccc}
Q=2M-2n_2L, & t_B=t_0+n_1l, & R=3M-3n_1L.
\end{array}
$$
\end{itemize}
\end{lemma}

\noindent{\sc Proof.} Suppose first $\eps>0$. Then $\eps\ge n_1$ because of
irreducibility of $t_B$. Since $R$ is irreducible and curves of the class
$t_0+n_2l$ sweep any effective divisor of the class $M-n_3L$, we have
$(t_0+n_2l)\circ R=3n_2-3\eps\ge 0$, i.e., $\eps\le n_2$.

Let $n_1=n_2=n_3$. Since $t_B\circ(M-n_1L)=0$, the linear system $|M-n_1L|$
is two-dimensional, and any of its elements is irreducible, there exists
at least one-dimensional subsystem in it that contains $t_B$.
Let $C\in|M-n_1L-t_B|$ be a general element, then
$$
    Q\circ C=S+S'\sim 2M^2-2\left(\frac13(b-\eps)+n_1\right)ML,
$$
where $S$ and $S'$ are ruled surfaces covering the base. Note that there
exists only a unique curve of the class $t_0$, and $t_B\cap t_0=\emptyset$.
At least one of these surfaces, say, $S$, containes both $t_0$ and $t_B$.
It is easy to see that $S$ is of the class $M^2-2n_1ML$, whence
$S'\sim M^2-\frac43n_1ML$. Since $S\ne S'$, then $t_0\not\subset S'$,
and denoting by $P$ the intersection of two general elements of
$|M-n_1L|$, we can see that $\dim(S'\cap P)=0$ (note that 
$\Bas|M-n_1L|=t_0$). We get a contradiction:
$$
 0\le S'\circ P=(M^2-\frac43n_1ML)\circ(M^2-2n_1ML)=-\frac{n_1}3 <0.
$$

So if $\eps\ge n_1$ then $n_1<n_3$; moreover, $n_1>0$ because of
irreducibility of $R_Q$ (otherwise, if $n_1=0$, there should be at least
1-dimensional family of lines of the class $t_0$ sweeping a surface, 
and $R_Q$ should contain it).

Let $G$ be an irreducible ruled surface that is the complete 
intersection of general elements of $|M-n_2L|$ and $|M-n_3L|$.
Let us show that $t_B\subset G$. Indeed, $G$ has the type ${\bf F_{n_1}}$
with $t_0$ as its exceptional section. Assume the converse, i.e., 
$t_B\not\subset G$. Then, for a general fiber $S\in|L|$, the line
$G|_S$ does not lie on $Q|_S$. So $G\not\subset Q$, and we may suppose
that $G\circ Q=t_0+C$, where $C$ is a curve with no $t_0$ as a component.
It is easy to compute that
$$
    C\sim t_0+2\left(b-n_2-n_3-\frac{b-\eps}3\right)l,
$$
and
$$
      2\left(b-n_2-n_3-\frac{b-\eps}3\right) \ge n_1,
$$
which, taking into account $\eps\le n_2$, gets a contradiction: 
$n_1>\frac43n_1$.

So $t_B\subset G$, hence $n_2<n_3$. In fact, the equality should led
$\dim|M-n_3L|$, thus $\Bas|M-n_3L|=t_0$, and we could always find $G$ 
such that $t_B\not\subset G$.

By $T$ we denote a unique effective element of $|M-n_3L|$. We show
that $T|_Q=2G$. Indeed, assuming the converse, i.e., $t|_Q=G+G'$,
$G\ne G'$, and taking into account that $t_B\subset G\cap G'$ and
$t_0\subset G$, we see $t_0\not\subset G'$. Then
$$
G'|_G=t_0+\mbox{\#$(t_0\cap t_B=\eps-n_1)$ fibers}.
$$
Choose a general $S\in|M-n_1L|$. It is clear that $S|_T$ is a ruled 
surface of the type ${\bf F_{n_2}}$ with the exceptional section $t_0$,
and $G\cap S|_T=t_0$. Let $C=G'\cap S|_T$ be an irreducible curve.
We have $C\circ t_0=\eps-n_1$, $Q\circ T\circ S=t_0+C$, whence
$$
    C\sim t_0+2\left(n_2-\frac{b-\eps}3\right)l.
$$
This yields $C\circ t_0=n_2-\frac23(b-\eps)=\eps-n_1$, and we get
a contradiction again:
$$
   \eps=n_1+n_2-2n_3<0.
$$
So we have
$$
T\circ Q=2M^2-2\left(n_3+\frac{b-\eps}3\right)ML=2G=2(M^2-(n_2+n_3)ML),
$$
thus $2n_2=n_1+n_3-\eps$. If $\eps$ were greater than $n_1$, then
$R_Q$ should contain $G$, because $G$ is covered by curves of the class
$t_0+n_1l$. This contradicts to irreducibility of $R_Q$. So $\eps=n_1$
and then $2n_2=n_3$. In order to show $n_2\ge 3n_1$, we note that
$S_Q=S\cap Q$ is a conic bundle without degenerations, i.e., a ruled
surface of the type ${\bf F_{n_2}}$. Then $R$ cuts off an effective curve 
on $S_Q$ that does not contain $t_0$, and $t_0$ itself. We get
$3n_3-3n_1\ge 5n_2$, which yields $n_2\ge 3n_1$.

It remains to show that $n_1$ is even. Let $C\sim t_0+n_1l$ be a 
general irreducible curve on $G$, and $C_V$ its inverse image on $V$.
Then $C_V$ covers $C$ with the ramification divisor of degree
$t_B\circ C=n_1$. Clearly, this degree must be even. This completes
the case $\eps>0$.

Now let $\eps=0$. We may assume $n_1<n_3$. Take $T\in|M-n_3L|$. It is 
easy to observe that $T=\Proj{\mc E_T}$, where
$$
   {\mc E}_T={\mc O}\oplus{\mc O}(n_1)\oplus{\mc O}(n_2).
$$
Since $T$ can be covered by curves of the class $t_0+n_2l$ and
$Q\sim 2M-\frac23bL$ is irreducible, then $n_2\ge\frac{b}3$,
i.e., $2n_2\ge n_1+n_3$, and $n_2>n_1$ since $n_3>n_1$.

Further, let $C\sim t_0+n_1l$ be irreducible. Since 
$b=n_1+n_2+n_3>3n_1$, then $C\circ Q=2(n_1-\frac{b}3)<0$, and
therefore $C\subset Q$. A ruled surface 
$G\sim (M-n_2L)\circ(M-n_3L)=M^2-(n_2+n_3)ML$ has the type 
${\bf F_{n_2}}$ and can be covered by curves of the same class as $C$,
so $G\subset Q$. Therefore, $T|_Q=G+G'$ for some $G'$.

Let $M_T=M|_T$ be the tautological divisor on $T$, $L_T=L|_T$.
Then $G$ has the class $M_T-n_2L_T$ in $T$, 
$Q|_T\sim 2M_T-\frac23bL_T$, whence $G'\sim M_T-(\frac23b-n_2)L_T$.
Note that $\frac23b-n_2>n_1$ since $n_1<n_2$.

But $T$ has only one divisor of the form $M_T-xL_T$ for $x>n_1$:
this is $M_T-n_2L_T$. So $\frac23b-n_2=n_2$, hence $n_1+n_3=2n_2$.

Now, let $C\subset Q$ be a general irreducible curve of the class
$t_0+n_1l$, $C_V$ its inverse image on $V$. Then $C_V$ is the double
cover of $C$ branched along a divisor of degree $R\circ C=3n_1$,
whence $n_1$ and $n_3$ are even. Lemma \ref{struct_lemma} is proved.

Such a construction of $V$ allows us to
deal completely with intersections of cycles on $V$. Denote $G_V$,
$F$, and $H$ the inverse images of $G\sim (M-n_2L)\circ(M-n_3L)$,
$L_Q=L|_Q$, and $M_Q=M|_Q$ on $V$ respectively. Then, let $2s_0$ be
the pull back of $t_0$, and $2f$ the pull back of a line in a fiber 
of $L_Q$.

\begin{lemma} $\overline{\bf NE}(V)={\bf NE}(V)=\RA_+[s_0]\oplus\RA_+[f]$,
and:
\begin{itemize}
\item[1)] if $\eps=0$, then $K_V=-G_V+(\frac12n_1-2)F$, $H=2G_V+n_3F$,
$K_V^2=s_0+(4-n_2)f$, $s_0\circ F=f\circ G_V=1$, $s_0\circ G_V=-\frac12n_3$,
$f\circ F=0$.
\item[2)] if $\eps=n_1>0$, then $K_V=-G_V-(\frac12n_1+2)F$, $H=2(G_V+n_2F)$,
$s_{б}=s_0+\frac12n_1f$, $K_V^2=s_0+(4+\frac32n_1-n_2)f$,
$s_0\circ F=f\circ G_V=1$, $s_0\circ G_V=-\frac12n_3=-n_2$, $f\circ F=0$.
\end{itemize}
\end{lemma}

\noindent{\sc Proof.} It can be easily checked by the following
way. First blow up $t_B$ on $Q$ with an exceptional divisor $S$;
then take the double cover branched along a smooth divisor composed
from the pre-image of $R_Q$ and $S$; after that, contract the pre-image
of $S$ onto $s_B$.

\begin{lemma}
\label{lem_1}
Let $\eps=0$ (i.e., $s_B=s_0$). If $|n(-K_V)-mF|$, $m>0$, has no fixed
components, then $n_3=2$.
\end{lemma}

\noindent{\sc Proof.} Let $D\in|n(-K_V)-mF|$. Since
$$
D|_{G_V}\sim ns_0-\left(m+n\left(\frac{n_1}2+n_3-n_2-2\right)\right)f
$$
is an effective curve, we obtain
$$
   0<\frac{m}n\le 2+n_2-n_3-\frac{n_1}2=2-\frac{n_3}2.
$$
It only remains to take into account that $n_3$ is positive and even.
The lemma is proved.

\begin{lemma}
\label{lem_2}
Let $\eps=n_1>0$ (so $s_B\ne s_0$). If $|n(-K_V)+mF|$ has no fixed 
components, then $m>0$.
\end{lemma}

\noindent{\sc Proof.} Let $D\in|n(-K_V)+mF|$. Since
$$
D|_{G_V}\sim ns_0+[m+n(n_1-n_2+2)]f
$$
is an effective curve, we have $2+n_1-n_2+\frac{m}n\ge 0$. Using 
$n_2\ge 3n_1$, we see that $2+n_1-n_2\le 2-2n_1<-2$, whence 
$m>0$. The lemma is proved.

\subsection{Results about rigidity.}

\begin{theorem}
\label{dp1_th}
Conjecture \ref{main_conj} is true for smooth Mori fibrations on del Pezzo
surfaces of degree 1 (over $\POn$).
\end{theorem}

In order to prove this theorem, we shall use the following proposition:

\begin{proposition}
\label{prop1}
Let $V/\POn$ be a smooth Mori fibration on del Pezzo surfaces of degree 1,
$V'/S'$ another Mori fibration, $\chi: V\dasharrow V'$ a birational map.
Suppose that
$$
  {\mc D'}=|n'(-K_{V'})+\mbox{pull back of an ample divisor from the base }|
$$
is a very ample linear system. Suppose aslo
$$
{\mc D}=\chi^{-1}_*{\mc D'}\subset|n(-K_V)+mF|
$$
for some $m\ge 0$. 

Then $n=n'$ (so $\chi$ is birational over the base), possibly except the 
case $\eps=n_1=0$, $n_2=1$, $n_3=2$.
\end{proposition}

\noindent{\sc Proof of the proposition.} By lemma \ref{easy_lemma},
$\mu(V',{\mc D'})=n'=\alpha(V',{\mc D'})$ and 
$\mu(V,{\mc D})=n=\alpha(V,{\mc D})$. If $n=n'$, the assertion follows from
proposition \ref{th_Corti}.

Suppose $n>n'$. Then a log pair $K_V+\frac1n{\mc D}$ is not canonical, 
i.e., the linear system ${\mc D}$ has a maximal singularity over some
point. But using \cite{Pukh1}, we may say a little more.

Namely, let $D_1$, $D_2$ be general elements of ${\mc D}$. We observe that
$$
   D_1\circ D_2= \left\{
\begin{array}{ll}
n^2s_0+((4-n_2)n^2+2mn)f, & \mbox{if $\eps=0$} \\
n^2s_0+((4+\frac32n_1-n_2)n^2+2mn)f, & \mbox{if $\eps=n_1>0$}
\end{array}
\right.
$$
Put down $D_1\circ D_2=Z^h+Z^v$, where $Z^h$ and $Z^v$ are horizontal
and vertical cycles respectively. Then there exist a geometric discrete
valuation $\nu$ centered at a point $B_0$ over some point $t\in \POn$, 
and a positive
number $e=\nu({\mc D})-n\delta$, where $\delta$ is the canonical 
multiplicity with respect to $\nu$, such that for the component 
$Z^v_t\subset F_t$ of the cycle $Z^v$ in the fiber $F_t$ over the 
point $t$ we have
$$
   \deg Z^v_t < \left\{
\begin{array}{ll}
(4-n_2)n^2+2n\frac{e}{\nu(F_t)}, & \mbox{if $\eps=0$} \\
(4+\frac32n_1-n_2)n^2+2n\frac{e}{\nu(F_t)}, & \mbox{if $\eps=n_1>0,$}
\end{array}
\right.
$$
where $\deg Z^v_t\eqdef Z^v_t\circ(-K_V)$ (this is so-called {\it the 
supermaximal singularity condition}, in the notions of 
\cite{Pukh1}). Besides, there exist positive numbers $\Sigma_0$,
$\Sigma_0'\le\Sigma_0$, and $\Sigma_1$ such that multiplicities
$m^h$ and $m^v$ of the cycles $Z^h$ and $Z^v_t$ at $B_0$ satisfy
$$
(\Sigma_0+\Sigma_1)(m^h\Sigma_0+m^v\Sigma_0')>(2n\Sigma_0+n\Sigma_1+e)^2.
$$
Note that $m^h\le\deg Z^h\eqdef Z^h\circ F=n^2$ and $m^v\le 2\deg Z^v_t$.
So, in the case $\eps=n_1>0$ we have
$$
(\Sigma_0+\Sigma_1)((9+3n_1-2n_2)n^2\Sigma_0+4ne)>
(2n\Sigma_0+n\Sigma_1+e)^2,
$$
hence
$$
   (2n_2-3n_1-5)\Sigma_0(\Sigma_0+\Sigma_1)+(n\Sigma_1-e)^2<0,
$$
i.e., $2n_2<5+3n_1$. Taking into account the condition $n_2\ge 3n_1$,
we obtain $n_1\le 1$, which is impossible because $n_1$ must be even.

Now let $\eps=0$. Then
$$
(\Sigma_0+\Sigma_1)((8-2n_2)n^2\Sigma_0+4ne)>(2n\Sigma_0+n\Sigma_1+e)^2,
$$
whence
$$
   (2n_2-5)\Sigma_0(\Sigma_0+\Sigma_1)+(n\Sigma_1-e)^2<0.
$$
So $n_2<\frac52$, so the only three cases are possible:
\begin{itemize}
\item[] $n_1=0$, $n_2=2$, $n_3=4$;
\item[] $n_1=2$, $n_2=2$, $n_3=2$;
\item[] $n_1=0$, $n_2=1$, $n_3=2$.
\end{itemize}
We show that the first two cases can not occur. Suppose first that there are 
no sections of the class $s_0$ through the point $B_0$.
Let $Z^h_1$ be the sum of all horizontal cycles through $B_0$. So
$Z_1^h\sim As_0+Bf$ for some $B\ge A$. Since $n_2=2$ in the considering
cases, we have
$$
\sum_{p\in\POn}\deg Z^v_p \le 2n^2+2mn-B \le 2n^2+2mn-A,
$$
which implies the following inquality for the supermaximal singularity:
$$
m^v\le 2\deg Z_t^v < 4n^2+4n\frac{e}{\nu(F_t)}-2A.
$$
But we know that $m^h\le A$, and then
\begin{equation}
\label{eq1}
(\Sigma_0+\Sigma_1)(A\Sigma_0+(4n^2-2A)\Sigma_0+4ne)>
     (2n\Sigma_0+n\Sigma_1+e)^2,
\end{equation}
so we get a contradiction:
$$
   A\Sigma_0(\Sigma_0+\Sigma_1)+(n\Sigma_1-e)^2<0.
$$
It remains to prove that there are no maximal singularities over points
that lie on sections of the class $s_0$.

Consider the case $n_1=n_2=n_3=2$. Note that $s_B\sim s_0$ is a unique
section of such a class in this case. If ${\mc D}$ has a maximal
singularity over the point $B_0\in s_B$, then for a general $D\in{\mc D}$
it holds $\nu_1=\mult_{B_0}D>n$, and for a general curve $C\sim f$ in the
fiber through $B_0$ we get $n=D\circ C\ge\nu_1>n$, which is impossible.

Now let $n_1=0$, $n_2=2$, $n_3=4$. Proving of this case consists of two
lemmas below and an observation that any curve of the class $s_0$ lies on
the divisor $G_V$. Recall that the valuation $\nu$ can be realized as 
a chain of blow-ups centered in nonsingular centers $B_0, B_1,\ldots$;
a general $D\in{\mc D}$ has multiplicities $\nu_1,\nu_2,\ldots$ in these
centers respectively.

\begin{lemma}
\label{spec_lemma}
Let ${\mc D}\subset|n(-K_V)+mF|$, $m\ge 0$, has a maximal singularity
over a point $B_0\in S$, where $S$ is a fiber of $V/\POn$. By 
$l\in|-K_S|$ we denote the curve through $B_0$. Then
\begin{itemize}
\item[a)] $S$ is smooth at $B_0$;
\item[b)] $B_0$ is a double point of $l$;
\item[c)] $\nu$ defines an infinitely near singularity, and if 
	   $B_1\cap S^1\ne\emptyset$ (upper indices denote the strict
	   transform on the correspondig floor of the chain of
	   blow-ups realizing $\nu$), then $B_1\cap l^1\ne\emptyset$ and
	   $\nu_1+\nu_2\le 3n$.
\end{itemize}
\end{lemma}

\noindent{\sc Proof.} If $S$ is singular at $B_0$, then a general curve
$C\in|2(-K_S)-B_0|$ is also singular at this point, so we get
a contradiction:
$$
   2n=D\circ C\ge 2\nu_1>\nu_1.
$$
Further, suppose $l$ is nonsingular at $B_0$. Set 
$$
   D|_S\sim kl+C,
$$
where $l\not\subset\Supp C$. We have
$$
    n < \nu_1 \le \mult_{B_0}D|_S=k+\mult_{B_0}C.
$$
But $\mult_{B_0}C\le C\circ l=(n(-K_S)-kl)\circ l=n-k$, and then
$$
        n < \nu_1\le k+n-k=n.
$$
So $l$ has a double point at $B_0$. Now we show that $\nu$ defines an
infinitely near singularity. Assuming the converse, we get $\nu_1>2n$.
Supposing $D|_S\sim kl+C$, we see that
$$
   2\mult_{B_0}C\le C\circ l=n-k,
$$
whence 
$$
 2n<\nu_1\le 2k+\mult_{B_0}C\le \frac32k+\frac12n,
$$
i.e., $k>n$, a contradiction.

Now let $B_1\cap S^1\ne\emptyset$, but $B_1\cap l^1=\emptyset$. Denote
$E_1$ the exceptional divisor of the blow-up of $B_0$. Then
$$
   D^1|_{S^1}=(D\circ S)^1+mE_1|_{S^1},
$$
and for the decomposition $D|_S\sim kl+C$ we get
$$
   \tilde\nu_1+\tilde\nu_2\le 2k+\mult_{B_0}C+\mult_{\tilde B_1}C^1+m,
$$
where $\tilde B_1=B_1\cap S^1$, $\tilde\nu_1=\mult_{B_0}D|_S$,
$\tilde\nu_2=\mult_{\tilde B_1}(D|_S)^1$. Then,
$$
  \mult_{\tilde B_1}C^1\le\mult_{B_0}C\le
      \frac12 (n(-K_S)-kl)\circ l=\frac12(n-k),
$$
so
$$
\tilde\nu_1+\tilde\nu_2\le 2k+2\frac{n-k}2+m=n+k+m.
$$
But $\tilde\nu_1\ge\nu_1+m$ and $\tilde\nu_2\ge\nu_2$, thus
$$
   2n < \nu_1+\nu_2\le n+k,
$$
i.e., $k>n$, which is impossible.

So $B_1\cap l^1\ne\emptyset$. In order to show that $\nu_1+\nu_2\le 3n$,
we argue as before, only taking into account that
$$
   \tilde\nu_1+\tilde\nu_2\le 3k+\mult_{B_0}C+\mult_{\tilde B_1}C^1+m.
$$
We have $\nu_1+\nu_2\le 2k+n\le 3n$. The lemma is proved.

\begin{lemma}
Let $\eps=0$, $n_1=0$. Then $G_V\cong C\times\POn$, where $C$ is an
elliptic curve.
\end{lemma}

\noindent{\sc Proof.} We use the notation of section \ref{subsec1}. It is
obvious that $G_V\cong C\times\POn$. We shall prove the smoothness of $C$.

Notice that $R|_G\sim 3t_0$. Let $B$ be an irreducible component
of $\Supp R|_G$. Clearly, $B\sim t_0$. Let $\psi:\tilde V\to V$ be
the blow-up of $B$ with the exceptional divisor $E$. $E$ is a ruled
surface of the type ${\bf F_{n_2}}$. Denote $t_E$ and $l_E$ the classes
of the exceptional section and a fiber of $E$ respectively. Then 
for strict transforms of the divisors $G$ and $R$ on $\tilde V$ we have
(using their smoothness along $B$)
$$
\begin{array}{cc}
    \tilde G|_E\sim t_E,  &  \tilde R|_E\sim t_E+n_2l_E,
\end{array}
$$
hence $\tilde G|_E\cap\tilde R|_E=\emptyset$. So any fiber of 
$G\cong\POn\times\POn$ meets $R$ transversally, i.e., at three different
points. Thus, $C$ is nonsingular.

The lemma is proved. This completes also the proof of proposition \ref{prop1}.

\begin{remark}
\label{rem_1}
In fact, proposition \ref{prop1} works also in the case 
$n_1=0$, $n_2=1$, $n_3=2$, but the proof is more complicated. Since this
case is non-rigid, satisfies theorem \ref{dp1_th}, and can not be advanced
to a more complete description yet (see below), we omit proving of it.
\end{remark}

\begin{corollary}
\label{dp1_cor}
For $\eps=n_1>0$, all smooth Mori fibrations on del Pezzo surfaces of 
degree 1 over $\POn$ are rigid over the base; for $\eps=0$, all cases
are also rigid, except the following:
\begin{itemize}
\item[a)] $n_1=n_2=n_3=2$;
\item[b)] $n_1=0$, $n_2=1$, $n_3=2$.
\end{itemize}
\end{corollary}

\noindent{\sc Proof.} Rigidity follows from proposition \ref{prop1}, using
lemma \ref{lem_1} in the case $\eps=0$ or lemma \ref{lem_2} when 
$\eps=n_1>0$. 

We will deal with non-rigid cases after the following remark.

\begin{remark}
In \cite{Pukh1} it was proved that (smooth) pencils of del Pezzo surfaces 
(of degree 1, 2 or 3) are all rigid, if the so-called 
${\bf K^2}$-condition holds:
\begin{center}
{\it
cycles $aK^2-bf$ are not effective for any $a,b>0$.
}
\end{center}
For degree 1, this condition is satisfied exactly when $n_2\ge 4$ (for
$\eps=0$) or $n_2\ge 4+\frac32n_1$ (for $\eps=n_1>0$). The second is always 
true except $n_1=2$, $n_2=6$, $n_3=12$.

Another sufficient condition of rigidity was proposed in \cite{Isk1} 
(conjecture 1.2): 
\begin{center}
{\it
$(-K_V)^3+m_0+1\le 2$},
\end{center}
where
$$
    m_0=\min\{r\in\QA: \mbox{$(-K_V)+rF$ is nef}\}.
$$
It holds exactly when $n_2\ge 3$ ($\eps=0$) or $n_2\ge 2+\frac32n_1$
($\eps=n_1>0$). Notice that the second is always true.

Thus, the second sufficient condition is more exact. Both these conditions
do not include the rigid case $n_1=0$, $n_2=2$, $n_3=4$ for $\eps=0$.
\end{remark}

Now we consider the non-rigid cases. We start with the case $n_1=0$,
$n_2=1$, $n_3=2$. It was already considered in some works (see 
\cite{Isk1}, section 2, \S 2, and \cite{Kha1}). It is easy to see that
a unique effective divisor of the class $G_V$ is the direct product of
$\POn$ and an elliptic curve, the linear system 
${\mc H}=|2G_V+2F|=|2(-K_V)-2F|$ is base points free and defines 
a morphism contracting $G_V$ along $\POn$ onto a curve $l$ on a Fano
variety $U$ of index 2 and degree 5 (this is so-called the double cone
over the Veronese surface). Indeed, the linear system $|M_Q|$ maps
$Q$ onto a cone $\tilde Q\in\PS$ over the Veronese surface in $\PF$, and
$U$ is obtained by taking the double cover branched along a cubic section
that does not pass through the cone vertex. $l$ covers one of the 
generators of $\tilde Q$. Let
$$
       \mu_l: V=V_l\to U
$$
be such a contraction. Notice that $-K_V\sim G_V+2F$, and $V$ is a Fano
itself. Thus, $V/\POn$ is not rigid, and 
$1=\alpha(V,{\mc H})<\mu(V,{\mc H})=2$, i.e., the conditions of 
conjecture \ref{main_conj} hold. Notice that $U$ has a lot of 
such structures. Unfortunately, the technique of the maximal singularities
method is not enough yet to deal completely with this case (for example,
it fails when $l$ has a double point; see also remark \ref{rem_1}).
Nevertheless, we can formulate the following conjecture, which seems to 
be true:

\begin{conjecture}
\label{conj_1}
Any smooth Mori fibration in the class of birational equivalency of $U$ is
biregular to either $U$ or $V_l$ for some $l$.
\end{conjecture}

The remaining case is $\eps=0$, $n_1=n_2=n_3=2$. Let $V/\POn$ be the 
corresponding pencil of del Pezzo surfaces of degree 1. Notice that
$$
    {\mc N}_{s_0|V}\simeq{\mc O}(-1)\oplus{\mc O}(-1)
$$
and $|mH|$  for $m\gg 0$ gives a contraction of $s_0$. So there exists
a flop
$$
        \mu: V\dasharrow U
$$
centered at $s_0$. Then, $\dim|G_V|=1$:
$$
 2\dim|G_V|=2\dim|G|\le\dim|2G|=\dim\left|(M-2L)|_Q\right|=2.
$$
Let $G_U=\overline{\mu(F)}$, and $F_U$ the class of a fiber of $U/\POn$.
Obviously, $\mu^{-1}|G_U|=|F|$ and $\mu^{-1}|F_U|=|G_V|$. Moreover,
$U/\POn$ is a pencil of del Pezzo surfaces of degree 1 with the same
structure parameters $n_1=n_2=n_3=2$.

\begin{proposition}
\label{prop_flop}
Let $\chi:V\dasharrow W$ be a birational map onto a Mori fibration 
$\gamma:W\to S$. Then either $\chi$ or $\chi\circ\mu^{-1}:U\dasharrow W$
are birational over the base. Moreover, $U$ and $V$ are unique smooth
Mori models in their class of birational equivalency; $Bir(V)=Aut(V)$,
and in general case, $Bir(V)=<\sigma>_2$, where $\sigma$ is the double
cover involution.
\end{proposition}

\noindent{\sc Proof.} Let  ${\mc D}_W=|n'(-K_W)+\gamma^*(A)|$ be a very
ample linear system, where $A$ is an ample divisor on $S$, and 
${\mc D}_V=\chi^{-1}_*{\mc D}_W\subset|n(-K_W)+mF|$.

Suppose $m\ge 0$. Then from proposition \ref{th_Corti} it follows that
$n=n'$ and $\chi$ is birational over the base. So we may assume
$m=-l<0$. Denote ${\mc D}_U=\mu^{-1}{\mc D}_V$. Since $-K_U=G_U+F_U$,
we have
$$
   {\mc D}_U\subset|(n-l)(-K_U)+lF_U|,
$$
and $\chi\circ\mu^{-1}$ is birational over the base by proposition 
\ref{prop1}.

Now let $\chi:V\dasharrow W$ be a birational map onto a Mori fibration 
$W/\POn$ with $-K_W$ to be {\it nef}. We may assume $\chi$ to be birational 
over the base. Suppose ${\mc Y}=|n(-K_W)+m'F_W|$ is a very ample linear
system, $m'>0$, and ${\mc D}=\chi^{-1}_*{\mc Y}\subset|n(-K_V)+mF|$ for
some $m\ge 0$. If $K_V+\frac1n{\mc D}$ were canonical, i.e., $D$ had no
maximal singularities, then $\chi$ would be an isomorphism by proposition
\ref{th_Corti}.

So let ${\mc D}$ have maximal singularities. We may assume that they are 
all infinitely near (see \cite{Pukh1}). By $Y$ and $D$ we denote general
elements of ${\mc Y}$ and ${\mc D}$ respectively. Choose a resolution of
singularities of $\chi$:
$$
\begin{CD}
@.@.Z @.@.\\
@.\stackrel{\phi}{}\swarrow@.@.\searrow {\stackrel{\psi}{}}\\
V @.@.\stackrel{\chi}{\dasharrow} @.\quad\; W
\end{CD}
$$
Notice that
$$
   nK_Z+\psi^{-1}Y=m'\psi^*(F_W)+\sum a_iE'_i
$$
and
$$
   nK_Z+\phi^{-1}D=m\phi^*(F)+\sum b_iE_i,
$$
where $E'_i,E_i$ are exceptional divisors, $a_i$ are all positive. 
Since ${\mc D}$ has maximal singularities, then
$$
    {\mc M}=\{i: b_i<0\}
$$
is not empty. For a point $t\in\POn$ we denote 
${\mc M}_t=\{i: \phi(E_i)\in F_t\}$.

Obviously, $\dim|nK_Z+\psi^{-1}Y|=m'$, so
$$
   \dim|m\phi^*(F)+\sum b_iE_i|=\dim|(m'+m-m')\phi^*(F)+\sum b_iE_i|=m'.
$$
Denote $I=\{t\in\POn: {\mc M}_t\ne\emptyset\}$. Then there exists a split
$$
    m-m'=\sum_{t\in I}k_t,
$$
where $k_t$ are all positive, such that 
$$
   \dim|\sum_{t\in I}k_t\phi^{-1}(F_t)+
     \sum_{t\in I}\sum_{i\in{\mc M}_t}(k_tc_i-b_i)E_i|=0,
$$
where for $t\in I$ positive numbers $c_i$ are defined by
$$
   \phi^*(F_t)=\phi^{-1}_*(F_t)+\sum_{i\in{\mc M}_t}c_iE_i.
$$
So for any $t\in I$ and $i\in {\mc M}_t$ we have $k_tc_i-b_i\ge 0$,
and, for every i, $k_i$ is the smallest positive number with such a
property. Thus
$$
   k_t=\max_{i\in{\mc M}_t}\ulcorner b_ic_i^{-1}\urcorner,
$$
where $\ulcorner\;\urcorner$ is the round up.

Let $D_1,D_2\in{\mc D}$ be general. We can put down
$$
D_1\circ D_2=Z^h+Z^v\sim n^2s+(2n^2+2mn)f,
$$
where $Z^v$ and $Z^h$ are vertical and horizontal cycles. Since $s_0$
is unique in its class and centers of maximal singularities can not lie
on $s_0$, then horizontal cycles being caught by at least one of these 
centers can put down as $As_0+Bf$ for some $B\ge A$. Taking into account that
$$
   m\le m'+\sum_{t\in I}k_t,
$$
we get an estimation
$$
  \deg Z^v=\sum_{t\in\POn}\deg Z^v_t\le 2n^2+2mn-A\le
       2n^2+2m'n+2n\sum_{t\in I}k_t-A
$$
(the degrees of vertical cycles are estimated by intersection with
$-K_V$). Thus there exist $p\in I$ and $j\in {\mc M}_p$ such that
$$
   \deg Z^v_p\le 2n^2+2m'n+2n\ulcorner b_jc_j^{-1}\urcorner-A.
$$
As in \cite{Pukh1}, \S 4, there exist positive numbers $\Sigma_0$,
$\Sigma_1$, and $\Sigma_0'$, where $\Sigma_0\ge\Sigma_0'$
and $c_j\ge\Sigma_0'$, such that the following inequality for 
multiplicities $m^h=\mult_{B_0}Z^h$ and $m^v=\mult_{B_0}Z^v_j$
of horizontal and vertical cycles at $B_0=\phi(E_j)$ holds:
\begin{equation}
\label{kv_eq}
(\Sigma_0+\Sigma_1)(\Sigma_0m^h+\Sigma_0'm^v)\ge
          (2n\Sigma_0+n\Sigma_1+b_j)^2.
\end{equation}
Notice that $\ulcorner b_jc_j^{-1}\urcorner\Sigma_0\le b_j+\Sigma_0$.
By lemma \ref{spec_lemma}, $\mc D$ has an infinitely near singularity
over the point $B_0\in F_p$. Consider two cases.

\noindent{\bf 1.} $B_1\cap F_p^1=\emptyset$, i.e., the center of the second
blow up in the chain realizing the maximal singularity, does not intersect
the strict transform of the fiber. In this case, since $B_1$ is a point,
we have $\Sigma_0'\le\frac12\Sigma_0$. Then
$$
   \Sigma_0'm^v\le 2\deg Z_p^v\frac{\Sigma_0}2\le
      (2n^2+2m'n+2n-A)\Sigma_0+2nb_j,
$$
and, using (\ref{kv_eq}) and an estimation $m^h\le A$, we get
$$
   n^2-2m'n-2n<0,
$$
i.e.,
\begin{equation}
\label{ner_m'1}
    \frac{m'}n > \frac12-\frac1n.
\end{equation}

\noindent{\bf 2.} Now let $B_1\cap S^1\ne\emptyset$. Here we can not state 
that $\Sigma_0'\le\frac12\Sigma_0$, but lemma \ref{spec_lemma} gives us
a good estimation $\nu_1+\nu_2\le 3n$. 

We know that $m^h\le A$ and
$$
   \Sigma_0'm^v\le 2\deg Z_p^v\Sigma_0\le
      (4n^2+4m'n+4n-2A)\Sigma_0+4nb_j.
$$
Substituting this in (\ref{kv_eq}), we obtain
$$
  (A-4m'n-4n)\Sigma_0(\Sigma_0+\Sigma_1)+(n\Sigma_1-b_j)^2<0.
$$
Notice that a condition $b_j>0$ is nothing but the Noether-Fano
inequality (\cite{Pukh2}): for some set of non-increasing numbers 
$\{r_i\}$, where $\sum r_i=\Sigma_0+\Sigma_1$, it holds
$$
   b_j=\sum r_i\nu_i-2n\Sigma_0-n\Sigma_1>0.
$$
Applying $\nu_1+\nu_2\le 3n$, we get
$$
   (n\Sigma_1-b_j)^2\ge \frac12n^2\Sigma_0(\Sigma_0+\Sigma_1).
$$
Thus
$$
   \frac12n^2+A-4m'n-4n<0,
$$
and since $A\ge 0$, than
\begin{equation}
\label{ner_m'2}
   \frac{m'}n>\frac18-\frac1n.
\end{equation}
Comparing (\ref{ner_m'1}) and (\ref{ner_m'2}), we may assume that 
(\ref{ner_m'2}) holds always. It only remains to use the condition that
$-K_W$ is {\it nef}: we can choose ${\mc Y}$ such that $n$ is big enough
but $\frac{m'}{n}$ is as small as we want, and then get a contradiction
because of (\ref{ner_m'2}).

Now, substituting $U$ for $W$, we obtain the statment about the birational
automorphisms group. The uniqueness of smooth models of $U$ and $V$ 
follows from the description of all smooth (rigid and non-rigid) Mori 
fibrations on del Pezzo surfaces of degree 1 given above. 

Proposition \ref{prop_flop} is proved. This also completes corollary
\ref{dp1_cor} and theorem \ref{dp1_th}.

\section{Smooth varieties with a pencil of del Pezzo surfaces of degree 2.}

\subsection{The essential construction.} 
In this section we study smooth Mori fibrations on del Pezzo surfaces of 
degree 2. Let $\rho:V\to\POn$ be such a fibration. Denoting $F$ the class 
of a fiber of $\rho$, we have
$$
    \Pic(V)=\ZA[-K_V]\oplus\ZA[F]
$$
and $(-K_V)^2\circ F=2$. Since $\rho$ is flat and $-K_V$ is $\rho$-ample,
for some integer $m$ we have
$$
    \rho_*{\mc O}(-K_V+mF)={\mc E},
$$
where ${\mc E}={\mc O}\oplus{\mc O}(n_1)\oplus{\mc O}(n_2)$ is a vector
bundle of rank 3 over $\POn$, $0\le n_1\le n_2$. Set
$$
    X=\Proj{\mc E}
$$
with a natural projection $\pi:X\to\POn$. By $M$ we denote the class of
the tautological bundle on $X$, $L$ the class of a fiber. Then
$$
   \Pic(X)=\ZA[M]\oplus\ZA[L].
$$
It is easy to see that there exists a double cover $\phi:V\to X$
branched along a smooth divisor $R$ of the class $4M+2aL$ such that
$$
\rho=\pi\circ\phi.
$$
Further, let $t_0$ be the class of a section that corresponds to
${\mc E}\to{\mc O}\to 0$, $l$ the class of a line in a fiber of $\pi$,
$s_0=\frac12\phi^*(t_0)$, and $f=\frac12\phi^*(l)$. Then
$$
   \overline{\bf NE}(X)={\bf NE}(X)=\RA_+[t_0]\oplus\RA_+[l],
$$
$$
   \overline{\bf NE}(V)={\bf NE}(V)=\RA_+[s_0]\oplus\RA_+[f].
$$
Denote $b=n_1+n_2$, $H=\phi^*(M)$; clearly, $F=\phi^*(L)$. We have
$M^3=b$, $M^2=t_0+bl$, $K_X=-3M+(b-2)L$, $M\circ t_0=L\circ l=0$,
$M\circ l=L\circ t_0=1$, и на $V$: $K_V=-H+(a+b-2)F$, $H^2=2s_0+2bf$,
$H\circ F=2f$, $H\circ s_0=F\circ f=0$, $H\circ f=F\circ s_0=1$,
$(-K_V)^2=2s_0+(8-4a-2b)f$, $H^3=2b$, $(-K_V)^3=12-6a-4b$ (see \cite{Isk1}).

Such varieties with $K_V^2=2s_0+\beta f$ for $\beta\le 0$ were first
studied in \cite{Pukh1}. That was an exclusively important step in
studying of geometry of Fano-fibered varieties.

\subsection{Results about rigidity.}

\begin{theorem}
\label{dp2_th}
Let $V/\POn$ be a smooth Mori fibration on del Pezzo surfaces of degree
2 over $\POn$ with $K_V^2=2s_0+\beta f$, where $\beta\le 2$ (i.e.,
$2a+b\ge 3$). Then $V/\POn$ is birationally rigid.
\end{theorem}

\noindent{\sc Proof.} Let $\chi:V\dasharrow V'$ be a birational map onto
a Mori-fibration $\rho':V'\to S'$, ${\mc D'}=|-n'K_{V'}+\rho'{}^*A'|$
a linear system as in proposition \ref{th_Corti}, and
${\mc D}=\chi^{-1}_*{\mc D'}\subset|-nK_V+mF|$. We will denote $D$ 
a general element of ${\mc D}$. Obviously, ${\mc D}$ has no fixed 
components. If $n=n'$, the assertion follows from proposition 
\ref{th_Corti}. So let $n>n'$.

\begin{lemma}
\label{dp2_lem}
$m\ge 0$, i.e. $\mu({\mc D})=\alpha({\mc D})$.
\end{lemma}

\noindent{\sc Proof of the lemma.}  The cases $\beta\le 0$ are proved in 
\cite{Pukh1} (the so-called ${\bf K^2}$-condition holds in these cases). 
Let $\beta=2$. Then $2a+b=3$, $b$ is odd, so $n_1<n_2$. Clearly,
$\mu({\mc D})=n$.

Assume the converse, i.e., $m<0$. Consider a general curve
$C\sim 2s_0+2n_1f$ on $V$ (there is at least 1-dimensional family of
such curves). Then for general $D$ we have $D\circ C\ge 0$, so
$$
  (nH-(\frac{n}{2}(b-1)-m)F)\circ (s_0+n_1f)\ge 0,
$$
and we get a contradiction:
$$
   0 < \frac{-m}{n}\le n_2-\frac{b}{2}+\frac12=\frac{n_1-n_2+1}{2}\le 0.
$$
The lemma is proved.

By the lemma, we may assume $m\ge 0$. 
The following argumentation proving the theorem is nearly the same as in
\cite{Pukh1}.

Obviously, ${\mc D}$ is not canonical, so it has maximal singularities.
Let a curve $B$ be the center of a maximal singularity. Then $B$ is a 
section of $\rho$ not lying on the ramification divisor, and we can
"untwist" such a maximal singularity using the {\it Bertini} involution
centered at $B$.

So we assume that the center of any maximal singularity is a point.
Set $D_1\circ D_2=Z^h+Z^v$, where $Z^h$ and $Z^v$ are horizontal and
vertical (effective) cycles. We see that
$$
  Z^h+Z^v=2s_0+(\beta n^2+4mn)f.
$$
Denote $\deg Z^h\eqdef Z^h\circ F$ and $\deg Z^v\eqdef Z^v\circ(-K_V)$.
Then, there exist a discrete valuation $\nu$ centered over a point
$B_0$ in a fiber $F_t$ (over a point $t\in\POn$), a positive number
$e=\nu({\mc D})-n\delta$, where $\delta$ is the canonical multiplicity
with respect to $\nu$, such that a component $Z^v_t$ of $Z^v$ lying in
$F_t$ has the degree
$$
   \deg Z^v_t < \left\{
\begin{array}{ll}
\beta n^2+4n\frac{e}{\nu(F_t)}, & \mbox{if $\beta > 0;$} \\
  4n\frac{e}{\nu(F_t)}, & \mbox{if $\beta\le 0.$}
\end{array}
\right.
$$
Notice that $m^h\eqdef\mult_{B_0}Z^h\le\deg Z^h=2n^2$ and
$m^v\eqdef\mult_{B_0}Z^v_t\le\deg Z^v_t$. Then, there exist positive
numbers $\Sigma_0\ge\Sigma_0'$ and $\Sigma_1$ suct that the following
inequality holds:
$$
(\Sigma_0+\Sigma_1)(m^h\Sigma_0+m^v\Sigma_0')>(2n\Sigma_0+n\Sigma_1+e)^2.
$$
It only remains to substitute the estimations of $m^h$ and $m^v$: 
$$
(\Sigma_0+\Sigma_1)((2+\beta n^2)\Sigma_0+4ne) >
          (2n\Sigma_0+n\Sigma_1+e)^2,
$$
i.e.,
$$
(2-\beta)n^2\Sigma_0(\Sigma_0+\Sigma_1)+(n\Sigma_1-e)^2<0,
$$
which is impossible, if $\beta\le 2$. Theorem \ref{dp2_th} is proved.

Now we shall consider smooth Mori fibrations on del Pezzo surfaces of
degree 2 for $2a+b\le 2$. First, let $2a+b=2$, i.e., $K_V^2=2s_0+4f$.

\begin{lemma}
\label{lem_17}
For $2a+b=2$, the only following cases may occur:
\begin{itemize}
\item[1)] $a=-1$, $n_1=n_2=2$;
\item[2)] $a=1$, $n_1=n_2=0$;
\item[3)] $a=0$, $n_1=n_2=1$;
\item[4)] $a=0$, $n_1=0$, $n_2=2$;
\item[5)] $a=-2$, $n_1=2$, $n_2=4$;
\item[6)] $a=-3$, $n_1=2$, $n_2=6$;
\item[7)] $a=-4$, $n_1=2$, $n_2=8$;
\end{itemize}
The first three cases are non-rigid.
\end{lemma}

\noindent{\sc Proof.} Since $b\ge 0$, the only cases 2), 3), and 4) are
possible if $a\ge 0$.

Let $a<0$. Then $b>0$, and since $R\circ t_0<0$, any curve of the class 
$s_0$ lies on $R$. So $n_1>0$ because of irreducibility of $R$, and such
a curve is unique on $X$. Further, let $\psi:\tilde X\to X$ be the blow-up
of $t_0$ with $E$ as an exceptional divisor. $E$ is a ruled surface
of the type ${\bf F_{n_2-n_1}}$. By $t_E$ and $l_E$ we denote the classes
of an exceptional section and a fiber of $E$. Suppose $\tilde R$ is the 
strict transform of $R$, then
$$
    \tilde R|_E\sim t_E+(2a+n_2)l_E
$$
is an irreducible curve on $V$ because of the smoothness of $R$. So either
$2a+n_2=0$, or $2a+n_2\ge n_2-n_1$. It is easy to check that the second 
case is possible only if $n_1=n_2=2$ and $a=-1$. Further, if $2a+n_2=0$, 
then $n_1=2$. Moreover, curves of the class $t_0+n_1l$ lie on a unique
effective divisor of the class $M-n_2L$, and since $R$ is irreducible,
we have $R\circ(t_0+n_1l)\ge 0$, whence $a\ge -4$.

Now we show that the first three cases are non-rigid. Let $a=1$, 
$n_1=n_2=0$. We see that $X\cong\PTw\times\POn$, and $V$ is a conic 
bundle with respect to the projection onto $\PTw$: double covers of 
curves of the class $t_0$ are conics since $R\circ t_0=2$. So
$V$ is non-rigid. Notice that this projection is given by a linear
system $|-K_V-F|$, which is free from base points.

Let $a=-1$, $n_1=n_2=2$. A unique curve of the class $t_0$ on $X$ lies
on the ramification divisor $R$, so $s_0$ is unique on $V$, too. Note
that
$$
{\mc N}_{s_0|V}={\mc O}(-1)\oplus{\mc O}(-2)
$$
and $|nH|$ for $n\ge 2$ gives a birational morphism contracting $s_0$.
Then, the base set of the pencil $|M-L|$ is exactly $s_0$; all elements
of this pencil are smooth and isomorphic to ${\bf F_2}$. For a general
$S\in|M-L|$ the restriction $R|_S$ is composed from $s_0$ and some
three-section that does not intersect $s_0$, so after taking the double
cover and contracting the pre-image of $s_0$ the surface $S$ becomes
a del Pezzo surface of degree 1. This means that the anti-flip
$V/\POn\dasharrow V'/\POn$ (centered at $s_0$) gives us a Mori fibration on
del Pezzo surfaces of degree 1 with a terminal singular point lying on
the exceptional curve of $V'$. Thus, $V/\POn$ is not rigid. Note that
the pencil $|-K_V-F|$ has no fixed components.

Finally, for $a=0$, $n_1=n_2=1$ the linear system $|-2K_V|$ gives 
a small contraction onto the canonical model of $V$, which can be realized 
as double covering of a non-degenerated quadratic cone in $\PQ$ branched
along a quartic section. If a curve of the class $s_0$ is unique on $V$
(this means that $t_0$ lies on the ramification divisor), then $s_0$ is
{\bf -2}-curve of the width 2 (in the notions of \cite{Reid}). 
Otherwise, there are two curves of the class $s_0$, which are
disjoint and {\bf -2}-curves of the width 1. In both the cases we 
obtain another structure of a smooth Mori fibration on del Pezzo surfaces 
of  degree 2 after making a flop centered at these curves. Notice again
that $|-K_V-F|$ has no fixed components. The case of two curves was
studied in detail in \cite{Grin}.

Lemma \ref{lem_17} is proved.

\begin{remark}
It is highly likely that cases 4) -- 7) are all rigid. Unfortunately,
the author can not prove it yet. As it often occur in the practice of
the maximal singularity method, the problems are related to excluding 
of infinitely near singularities over points on some rational curves
of a special kind.
\end{remark}

Now we consider cases when $2a+b=1$, i.e., $K_V^2=2s_0+6f$.

\begin{lemma} 
If $2a+b=1$, the only following three cases may occur:
\begin{itemize}
\item[1)] $a=0$, $n_1=0$, $n_2=1$;
\item[2)] $a=-1$, $n_1=1$, $n_2=2$;
\item[3)] $a=-2$, $n_1=1$, $n_2=4$.
\end{itemize}
The first two are non-rigid.
\end{lemma}

\noindent{\sc Proof.} Arguing as before, it is easy to show that
$a\ge -2$, so the only 1) -- 3) are possible.

In order to show that the case 1) is non-rigid, note that $V$ can
be obtained by blowing up of an elliptic curve of degree 2 on a smooth
Fano variety $U$ of genus 9 and index 2 (this is so-called {\it the double
space of index 2}, i.e., the double cover of $\PT$ branched along a quartic).
Observe that the birational morphism $V\to U$ is defined by the linear
system $|-2K_V-2F|$. Some results (but very incomplete) about these 
varieties are contained in \cite{Kha2}.

Consider the case $a=-1$, $n_1=1$, $n_2=2$. There exists a unique curve 
of the class $s_0$, and it is a {\bf -2}-curve of the width 1. Let
$V\dasharrow V^+$ be a flop centered at this curve. A linear
system $|H-2F|$ contains a unique element, which we denote $G_V$.
Its strict transform $G_V^+$ on $V^+$ is  a surface that is isomorphic
to the double cover of $\PTw$ branched along either a smooth conic or
a couple of different lines. So $G_V^+$ is either $\POn\times\POn$, or
a quadratic cone in $\PT$. Then, there exists an extremal contraction
$V^+\to U$ of $G_V^+$. U is a double cone over the Veronese surface, but
with a quadratic singularity (arising from a quadratic singularity of 
the ramification divisor). We see also that $U$ has (birationally) 
structures of fibrations on del Pezzo surfaces of degree 1 (see section
\ref{sec2}, the case $n_1=0$, $n_2=1$, $n_3=2$). Notice that $|-K_V-F|$
has no fixed components.

The lemma is proved.

As after lemma \ref{lem_17}, the author can say the same thing about the
case 3): it should be rigid, but I can not prove it yet.

We complete this survey of del Pezzo fibrations by the following lemma:

\begin{lemma}
Cases $2a+b\le 0$ can not occur.
\end{lemma}

\noindent{\sc Proof.} The same reasons as above, except the case $a=0$.
If $a=0$, then $b=0$, so $V$ is isomorphic to the direct product of
$\POn$ and a smooth del Pezzo surface of degree 2. But in this case 
$V$ is not a Mori fibration because of the relative Picard number 
(it is equal to 8). The lemma is proved.

\end{document}